\documentclass[12pt]{amsart}

\makeatletter
\newcommand\mathcircled[1]{%
	\mathpalette\@mathcircled{#1}%
}
\newcommand\@mathcircled[2]{%
	\tikz[baseline=(math.base)] \node[draw,ellipse,inner sep=1pt] (math) {$\m@th#1#2$};%
}

\usepackage{multicol}
\usepackage{latexsym}
\usepackage{psfrag}
\usepackage{amsmath}
\usepackage{amssymb}
\usepackage{epsfig}
\usepackage{amsfonts}
\usepackage{amscd}
\usepackage{mathrsfs}
\usepackage{graphicx}
\usepackage{enumerate}
\usepackage[autostyle=false, style=english]{csquotes}
\MakeOuterQuote{"}
\usepackage{ragged2e}
\usepackage[all]{xy}
\usepackage{mathtools}
\usepackage{comment}

\theoremstyle{remark}

\newlength\ubwidth

\usepackage[dvipsnames]{xcolor}

\usepackage{placeins}
\usepackage{tikz}
\usetikzlibrary{shapes, positioning}
\usetikzlibrary{matrix,shapes.geometric,calc,backgrounds}

\newtheorem{theorem}{Theorem}

\newtheorem{ex}[theorem]{Example}

\newtheorem{cor}[theorem]{Corollary}

\usepackage{amsmath,amssymb,amsthm}
\usepackage[dvipsnames]{xcolor}
\usepackage[colorlinks=true,citecolor=RoyalBlue,linkcolor=red,breaklinks=true]{hyperref}

\usepackage{graphicx}
\usepackage{mathdots} 
\usepackage[margin=2.5cm]{geometry}
\usepackage{ytableau}

\usepackage{cancel}
\usepackage{blkarray}
\usepackage{booktabs}
\usepackage{microtype}
\usetikzlibrary{arrows.meta}

\begin{document}

	\title[]
	{On Partitions with Constant-Sum Gap Sequences}

	\author[\"{O}mr\"{u}uzun Seyrek]{Hal\.{ı}me \"{O}mr\"{u}uzun Seyrek}
	\address{Hal\.{ı}me \"{O}mr\"{u}uzun Seyrek, Faculty of Applied Sciences, 
		Biruni University, Zeytinburnu, Istanbul, 34015, Turkey}
	\email{hseyrek@biruni.edu.tr}
	
	\subjclass[2010]{05A17, 05A15, 11P84}
	
\keywords{constant-sum gap sequences, partition generating function, $q$-series}
	
	\date{2026}
	
\begin{abstract}
	We introduce and study \emph{constant-sum gap partitions}: for a fixed
	nonnegative integer $C$, these are partitions
	$\lambda=(\lambda_1,\ldots,\lambda_k)$ whose successive differences
	$d_i=\lambda_{i+1}-\lambda_i$ satisfy $d_i+d_{k-i}=C$ for all valid
	$i$. We derive the generating function for this family, treating
	partitions of odd and even length separately, and compute it
	explicitly for $C=0, 1,2,3,4$. We then construct explicit bijections
	between constant-sum gap partitions with $2m+1$, respectively $2m$,
	parts and partitions into odd, respectively even, parts in which parts
	below the largest part have bounded multiplicity, together with an exact
	weight relation between a partition and its image. Finally, we give a
	geometric characterization of constant-sum gap partitions in terms of
	their Ferrers diagrams: removing a prescribed, linearly growing number
	of cells from each row renders the diagram symmetric about its middle
	row, for both the odd- and even-length cases.
\end{abstract}

\maketitle

	\section{Introduction}
	\label{secIntro}

An (ordinary) integer partition $\lambda$ is an expression of a 
non-negative integer by an unordered sum of positive integers~\cite{Andrews1998}.  
\begin{align*}
	n = \lambda_1 + \lambda_2 + \cdots + \lambda_m
\end{align*}
The parts $\lambda_i$'s may be written 
in non-increasing or non-decreasing order, 
and a limited number of zeros may be allowed 
to appear as parts.  
The number $n$ being partitioned is called the \emph{weight}, 
denoted by $\left\vert \lambda \right\vert$
and the number of parts $m$ is called the \emph{length} of the partition, denoted by $\ell(\lambda)$.  
For example, 4 has the following five partitions. 
\begin{align*}
	4, \quad
	1+3, \quad 
	2+2, \quad 
	1+1+2, \quad 
	1+1+1+1.
\end{align*}

Throughout this paper, we write partitions by separating their parts with commas rather than plus signs. Thus, for example, the partition $1, 2, 5, 5, 7$ denotes the partition $1+2+5+5+7$.

The study of integer partitions subject to conditions on the differences
between their parts has a long history, and continues to be a subject of great interest. We refer the interested reader to the following works: \cite{And69}, \cite{BK16}, \cite{Cha16}, \cite{CY18}, \cite{Schur1926}.

In this paper, we impose a condition on the entire sequence of \emph{successive}
differences (the gap sequence) of a partition. Specifically, we study
partitions $\lambda=(\lambda_1,\ldots,\lambda_k)$ whose gap sequence
$d_i=\lambda_{i+1}-\lambda_i$ satisfies $d_i+d_{k-i}=C$ for a fixed
nonnegative integer $C$, that is, gaps equidistant from the two ends of
the sequence sum to a constant. We call such partitions
\emph{constant-sum gap partitions} and denote the corresponding family by
$\mathcal P_C$. This condition is reminiscent of the equal-gap condition $d_i = d_{k+1-i}$ studied in our companion paper \cite{OS26}, where we investigated partitions whose successive differences form a palindromic sequence; the two families are, however, structurally distinct. Here we relax this exact
symmetry requirement to a constant-sum condition.

We first derive the generating function for $\mathcal P_C$, again split
according to the parity of the partition length (Theorem~\ref{thm:GF}),
and work out the resulting $q$-series explicitly for $C=0, 1,2,3,4$
(Corollaries \ref{C=0}, \ref{C=1}, \ref{C=2}, \ref{C=3}, \ref{C=4}). Our main results give explicit
bijections between $\mathcal P_C$, restricted to partitions of a fixed
odd or even length, and partitions into odd or even parts with bounded
multiplicity (Theorems~\ref{thm:odd-bijection} and~\ref{thm:even-bijection}
below), together with an explicit weight relation between a partition
and its image. We further give a geometric characterization of $\mathcal
P_C$ in terms of the Ferrers diagram: removing a prescribed, linearly
growing number of cells from each row renders the diagram symmetric
about its middle row (Theorems~\ref{thm:ferrers-odd} and~\ref{thm:ferrers-even}).
For general background on partition bijections and the techniques used
throughout, we refer the reader to Pak's survey~\cite{Pak06}.

\section{Generating Functions for Constant-sum Gap Partitions}

Let
	\[
	\lambda=(\lambda_1,\lambda_2,\ldots,\lambda_k),
	\qquad
	1\le \lambda_1\le \lambda_2\le \cdots\le \lambda_k,
	\]
	be an integer partition. The \emph{gap sequence} of $\lambda$ is
	\[
	d_i=\lambda_{i+1}-\lambda_i,
	\qquad
	1\le i\le k-1.
	\]

	Let $C\ge0$ be fixed. We say that $\lambda$ has a
	\emph{constant-sum gap sequence} if
	\[
	d_i+d_{k-i}=C,
	\qquad
	1\le i\le
	\left\lfloor\frac{k}{2}\right\rfloor.
	\]
	The set of all such partitions is denoted by $\mathcal{P}_C$.

\begin{theorem} \label{thm:GF}
	Let $\mathcal{P}_C^{\mathrm{odd}}$ and
	$\mathcal{P}_C^{\mathrm{even}}$ denote the subsets of
	$\mathcal{P}_C$ consisting of partitions with an odd and an even
	number of parts, respectively.
	
	Then
	\[
	\sum_{\lambda\in\mathcal{P}_C^{\mathrm{odd}}}q^{|\lambda|}
	=
	\sum_{m\ge0}
	\frac{
		q^{\frac{Cm^2+(C+4)m+2}{2}}
	}{
		1-q^{2m+1}
	}
	\prod_{j=1}^{m}
	\left(
	\sum_{r=0}^{C}q^{r(2j-1)}
	\right).
	\]
	
	Moreover, if $C$ is odd, then
	\[
	\mathcal{P}_C^{\mathrm{even}}=\varnothing.
	\]
	
	If $C$ is even, then
	\[
	\sum_{\lambda\in\mathcal{P}_C^{\mathrm{even}}}q^{|\lambda|}
	=
	\sum_{m\ge0}
	\frac{
		q^{\frac{(m+1)(C(m+1)+4)}{2}}
	}{
		1-q^{2m+2}
	}
	\prod_{j=1}^{m}
	\left(
	\sum_{r=0}^{C}q^{2rj}
	\right).
	\]
\end{theorem}

\begin{proof}
	We treat the odd- and even-length cases separately.
	
	First, let
	\[
	\lambda=(\lambda_1,\lambda_2,\ldots,\lambda_{2m+1})
	\in\mathcal{P}_C^{\mathrm{odd}},
	\]
	and put
	\[
	a=\lambda_1.
	\]
	Since $\lambda$ is a partition into positive parts, we have $a\geq1$.
	Let
	\[
	d_i=\lambda_{i+1}-\lambda_i,
	\qquad 1\leq i\leq 2m.
	\]
	The constant-sum condition gives
	\[
	d_i+d_{2m+1-i}=C,
	\qquad 1\leq i\leq m.
	\]
	Hence, $d_1,d_2,\ldots,d_m$ may be chosen freely from
	$\{0,1,\ldots,C\}$, while the remaining gaps are uniquely determined by
	\[
	d_{2m+1-i}=C-d_i,
	\qquad 1\leq i\leq m.
	\]
	
	Since
	\[
	\lambda_j
	=
	a+\sum_{i=1}^{j-1}d_i,
	\qquad 1\leq j\leq 2m+1,
	\]
	the weight of $\lambda$ is
	\[
	|\lambda|
	=
	(2m+1)a+\sum_{i=1}^{2m}(2m+1-i)d_i.
	\]
	Pairing the contributions of $d_i$ and $d_{2m+1-i}$, we obtain
	\[
	\begin{aligned}
		&(2m+1-i)d_i+i\,d_{2m+1-i}\\
		&\qquad
		=(2m+1-i)d_i+i(C-d_i)\\
		&\qquad
		=iC+(2m+1-2i)d_i.
	\end{aligned}
	\]
	Therefore,
	\[
	|\lambda|
	=
	(2m+1)a
	+
	C\sum_{i=1}^{m}i
	+
	\sum_{i=1}^{m}(2m+1-2i)d_i,
	\]
	and hence
	\[
	|\lambda|
	=
	(2m+1)a
	+
	C\binom{m+1}{2}
	+
	\sum_{i=1}^{m}(2m+1-2i)d_i.
	\]
	
	Set
	\[
	e_j=d_{m+1-j},
	\qquad 1\leq j\leq m.
	\]
	Then $0\leq e_j\leq C$, and
	\[
	|\lambda|
	=
	(2m+1)a
	+
	C\binom{m+1}{2}
	+
	\sum_{j=1}^{m}(2j-1)e_j.
	\]
	Conversely, every choice of
	\[
	a\geq1
	\qquad\text{and}\qquad
	e_1,e_2,\ldots,e_m\in\{0,1,\ldots,C\}
	\]
	determines uniquely the first $m$ gaps and hence, through the
	constant-sum condition, all the remaining gaps. It therefore determines
	a unique member of $\mathcal{P}_C^{\mathrm{odd}}$ with $2m+1$ parts.
	
	Thus, the generating function for the members of $\mathcal{P}_C$
	having exactly $2m+1$ parts is
	\[
	\begin{aligned}
		\sum_{a\geq1}
		q^{(2m+1)a+C\binom{m+1}{2}}
		\prod_{j=1}^{m}
		\left(
		\sum_{r=0}^{C}q^{r(2j-1)}
		\right)
		=
		\frac{
			q^{2m+1+C\binom{m+1}{2}}
		}{
			1-q^{2m+1}
		}
		\prod_{j=1}^{m}
		\left(
		\sum_{r=0}^{C}q^{r(2j-1)}
		\right).
	\end{aligned}
	\]
	Since
	\[
	2m+1+C\binom{m+1}{2}
	=
	\frac{Cm^2+(C+4)m+2}{2},
	\]
	summing over all $m\geq0$ gives
	\[
	\sum_{\lambda\in\mathcal{P}_C^{\mathrm{odd}}}q^{|\lambda|}
	=
	\sum_{m\geq0}
	\frac{
		q^{\frac{Cm^2+(C+4)m+2}{2}}
	}{
		1-q^{2m+1}
	}
	\prod_{j=1}^{m}
	\left(
	\sum_{r=0}^{C}q^{r(2j-1)}
	\right).
	\]
	
	We now consider the even-length case. Let
	\[
	\lambda=(\lambda_1,\lambda_2,\ldots,\lambda_{2m+2})
	\in\mathcal{P}_C^{\mathrm{even}},
	\]
	and again put
	\[
	a=\lambda_1\geq1.
	\]
	The gap sequence has length $2m+1$, and its middle gap is $d_{m+1}$.
	The constant-sum condition at the middle position yields
	\[
	d_{m+1}+d_{m+1}=C,
	\]
	or equivalently,
	\[
	2d_{m+1}=C.
	\]
	Hence, if $C$ is odd, no such partition exists, and therefore
	\[
	\mathcal{P}_C^{\mathrm{even}}=\varnothing.
	\]
	
	Suppose now that $C$ is even. Then
	\[
	d_{m+1}=\frac{C}{2}.
	\]
	The gaps $d_1,d_2,\ldots,d_m$ may be chosen freely from
	$\{0,1,\ldots,C\}$, and the remaining gaps are determined by
	\[
	d_{2m+2-i}=C-d_i,
	\qquad 1\leq i\leq m.
	\]
	As before,
	\[
	|\lambda|
	=
	(2m+2)a+
	\sum_{i=1}^{2m+1}(2m+2-i)d_i.
	\]
	For $1\leq i\leq m$, the paired gaps contribute
	\[
	\begin{aligned}
		&(2m+2-i)d_i+i\,d_{2m+2-i}\\
		&\qquad
		=(2m+2-i)d_i+i(C-d_i)\\
		&\qquad
		=iC+(2m+2-2i)d_i.
	\end{aligned}
	\]
	The middle gap contributes
	\[
	(m+1)d_{m+1}
	=
	\frac{C(m+1)}{2}.
	\]
	Therefore,
	\[
	\begin{aligned}
		|\lambda|
		={}
		(2m+2)a
		+
		C\sum_{i=1}^{m}i
		+
		\frac{C(m+1)}{2}
		+
		\sum_{i=1}^{m}(2m+2-2i)d_i.
	\end{aligned}
	\]
	Since
	\[
	C\sum_{i=1}^{m}i+\frac{C(m+1)}{2}
	=
	\frac{C(m+1)^2}{2},
	\]
	we obtain
	\[
	|\lambda|
	=
	(2m+2)a
	+
	\frac{C(m+1)^2}{2}
	+
	\sum_{i=1}^{m}2(m+1-i)d_i.
	\]
	
	Set
	\[
	e_j=d_{m+1-j},
	\qquad 1\leq j\leq m.
	\]
	Then $0\leq e_j\leq C$, and
	\[
	|\lambda|
	=
	(2m+2)a
	+
	\frac{C(m+1)^2}{2}
	+
	\sum_{j=1}^{m}2j e_j.
	\]
	Conversely, every choice of
	\[
	a\geq1
	\qquad\text{and}\qquad
	e_1,e_2,\ldots,e_m\in\{0,1,\ldots,C\}
	\]
	determines uniquely a partition in
	$\mathcal{P}_C^{\mathrm{even}}$ having $2m+2$ parts.
	
	Consequently, for even $C$, the generating function for the members
	having exactly $2m+2$ parts is
	\[
	\begin{aligned}
		\sum_{a\geq1}
		q^{(2m+2)a+\frac{C(m+1)^2}{2}}
		\prod_{j=1}^{m}
		\left(
		\sum_{r=0}^{C}q^{2rj}
		\right)=
		\frac{
			q^{2m+2+\frac{C(m+1)^2}{2}}
		}{
			1-q^{2m+2}
		}
		\prod_{j=1}^{m}
		\left(
		\sum_{r=0}^{C}q^{2rj}
		\right).
	\end{aligned}
	\]
	Finally,
	\[
	2m+2+\frac{C(m+1)^2}{2}
	=
	\frac{(m+1)(C(m+1)+4)}{2}.
	\]
	Summing over all $m\geq0$ gives
	\[
	\sum_{\lambda\in\mathcal{P}_C^{\mathrm{even}}}q^{|\lambda|}
	=
	\sum_{m\geq0}
	\frac{
		q^{\frac{(m+1)(C(m+1)+4)}{2}}
	}{
		1-q^{2m+2}
	}
	\prod_{j=1}^{m}
	\left(
	\sum_{r=0}^{C}q^{2rj}
	\right).
	\]
	This completes the proof.
\end{proof}

\begin{cor}
	The generating function for the family $\mathcal{P}_C$ is given by
	\[
	\sum_{\lambda\in\mathcal{P}_C}q^{|\lambda|}
	=
	\sum_{\lambda\in\mathcal{P}_C^{\mathrm{odd}}}q^{|\lambda|}
	+
	\sum_{\lambda\in\mathcal{P}_C^{\mathrm{even}}}q^{|\lambda|},
	\]
	where the second summand is interpreted as zero whenever $C$ is odd.
\end{cor}

\begin{proof}
	Since
	\[
	\mathcal{P}_C
	=
	\mathcal{P}_C^{\mathrm{odd}}
	\sqcup
	\mathcal{P}_C^{\mathrm{even}},
	\]
	the generating function of $\mathcal{P}_C$ is obtained by adding the
	generating functions of $\mathcal{P}_C^{\mathrm{odd}}$ and
	$\mathcal{P}_C^{\mathrm{even}}$. The result now follows from
	Theorem~\ref{thm:GF}.
\end{proof}

\begin{cor} \label{C=0}
	For $C=0$, the family $\mathcal{P}_0$ consists precisely of the
	partitions whose parts are all equal. Its generating function is
	\[
	\sum_{\lambda\in\mathcal{P}_0}q^{|\lambda|}
	=
	\sum_{k\geq1}\frac{q^k}{1-q^k}.
	\]
	Equivalently,
	\[
	\sum_{\lambda\in\mathcal{P}_0}q^{|\lambda|}
	=
	\sum_{n\geq1}\tau(n)q^n,
	\]
	where $\tau(n)$ denotes the number of positive divisors of $n$.
\end{cor} 

\begin{proof}
	Let
	\[
	\lambda=(\lambda_1,\ldots,\lambda_k)\in\mathcal{P}_0,
	\]
	with gap sequence
	\[
	d_i=\lambda_{i+1}-\lambda_i.
	\]
	Since
	\[
	d_i+d_{k-i}=0
	\]
	and every gap is nonnegative, each gap occurring in the condition must
	be zero. Hence,
	\[
	d_1=d_2=\cdots=d_{k-1}=0,
	\]
	so that
	\[
	\lambda=(a,a,\ldots,a)
	\]
	for some integers $a\geq1$ and $k\geq1$.
	
	Conversely, every partition of the form
	\[
	(a^k)=(a,a,\ldots,a)
	\]
	has zero gap sequence and therefore belongs to $\mathcal{P}_0$.
	Thus,
	\[
	\sum_{\lambda\in\mathcal{P}_0}q^{|\lambda|}
	=
	\sum_{k\geq1}\sum_{a\geq1}q^{ak}
	=
	\sum_{k\geq1}\frac{q^k}{1-q^k}.
	\]
	Finally, the coefficient of $q^n$ counts the pairs $(a,k)$ satisfying
	$ak=n$, and hence equals the number $\tau(n)$ of positive divisors of
	$n$.
\end{proof}

\begin{cor} \label{C=1}
	For $C=1$, the generating function for partitions in $\mathcal{P}_1$
	is
	\[
	\sum_{\lambda\in\mathcal{P}_1}q^{|\lambda|}
	=
	\sum_{m\geq0}
	\frac{
		q^{\frac{m^2+5m+2}{2}}
	}{
		1-q^{2m+1}
	}
	\prod_{j=1}^{m}
	\left(1+q^{2j-1}\right).
	\]
	Equivalently,
	\[
	\sum_{\lambda\in\mathcal{P}_1}q^{|\lambda|}
	=
	\sum_{m\geq0}
	\frac{
		q^{\frac{m^2+5m+2}{2}}
		(-q;q^2)_m
	}{
		1-q^{2m+1}
	}.
	\]
\end{cor}

\begin{proof}
	Since $C=1$ is odd, Theorem~\ref{thm:GF} implies that
	\[
	\mathcal{P}_1^{\mathrm{even}}=\varnothing.
	\]
	Substituting $C=1$ into the generating function for
	$\mathcal{P}_C^{\mathrm{odd}}$ gives
	\[
	\sum_{\lambda\in\mathcal{P}_1}q^{|\lambda|}
	=
	\sum_{m\geq0}
	\frac{
		q^{\frac{m^2+5m+2}{2}}
	}{
		1-q^{2m+1}
	}
	\prod_{j=1}^{m}
	\left(
	\sum_{r=0}^{1}q^{r(2j-1)}
	\right).
	\]
	Since
	\[
	\sum_{r=0}^{1}q^{r(2j-1)}
	=
	1+q^{2j-1},
	\]
	the first formula follows. The second follows from
	\[
	(-q;q^2)_m
	=
	\prod_{j=1}^{m}\left(1+q^{2j-1}\right).
	\]
\end{proof}

\begin{cor} \label{C=2}
	For $C=2$, the generating function for partitions in $\mathcal{P}_2$
	is
	\begin{align*}
	\sum_{\lambda\in\mathcal{P}_2}q^{|\lambda|}
	=
	\sum_{m\geq0}
	\frac{
		q^{m^2+3m+1}
	}{
		1-q^{2m+1}
	}
	\prod_{j=1}^{m}
	\left(
	1+q^{2j-1}+q^{2(2j-1)}
	\right)
	+\sum_{m\geq0}
	\frac{
		q^{(m+1)(m+3)}
	}{
		1-q^{2m+2}
	}
	\prod_{j=1}^{m}
	\left(
	1+q^{2j}+q^{4j}
	\right).
\end{align*}
	Equivalently,
	\[
	\sum_{\lambda\in\mathcal{P}_2}q^{|\lambda|}
	=
	\sum_{m\geq0}
	q^{m^2+3m+1}
	\frac{(q^3;q^6)_m}{(q;q^2)_{m+1}}
	+
	\sum_{m\geq0}
	q^{(m+1)(m+3)}
	\frac{(q^6;q^6)_m}{(q^2;q^2)_{m+1}}.
	\]
\end{cor}

\begin{proof}
	Since $C=2$ is even, both the odd- and even-length cases contribute.
	Substituting $C=2$ into Theorem~\ref{thm:GF} gives
	\[
	\sum_{\lambda\in\mathcal{P}_2^{\mathrm{odd}}}q^{|\lambda|}
	=
	\sum_{m\geq0}
	\frac{
		q^{m^2+3m+1}
	}{
		1-q^{2m+1}
	}
	\prod_{j=1}^{m}
	\left(
	\sum_{r=0}^{2}q^{r(2j-1)}
	\right),
	\]
	and
	\[
	\sum_{\lambda\in\mathcal{P}_2^{\mathrm{even}}}q^{|\lambda|}
	=
	\sum_{m\geq0}
	\frac{
		q^{(m+1)(m+3)}
	}{
		1-q^{2m+2}
	}
	\prod_{j=1}^{m}
	\left(
	\sum_{r=0}^{2}q^{2rj}
	\right).
	\]
	Since
	\[
	\sum_{r=0}^{2}q^{r(2j-1)}
	=
	1+q^{2j-1}+q^{2(2j-1)}
	\]
	and
	\[
	\sum_{r=0}^{2}q^{2rj}
	=
	1+q^{2j}+q^{4j},
	\]
	the first formula follows by adding the odd- and even-length
	generating functions.
	
	For the equivalent form, we use
	\[
	1+x+x^2=\frac{1-x^3}{1-x}.
	\]
	Thus,
	\[
	\prod_{j=1}^{m}
	\left(1+q^{2j-1}+q^{2(2j-1)}\right)
	=
	\frac{(q^3;q^6)_m}{(q;q^2)_m},
	\]
	and
	\[
	\prod_{j=1}^{m}
	\left(1+q^{2j}+q^{4j}\right)
	=
	\frac{(q^6;q^6)_m}{(q^2;q^2)_m}.
	\]
	Combining these products with
	\[
	\frac{1}{1-q^{2m+1}}
	\qquad\text{and}\qquad
	\frac{1}{1-q^{2m+2}},
	\]
	respectively, yields
	\[
	\frac{(q^3;q^6)_m}{(q;q^2)_{m+1}}
	\]
	and
	\[
	\frac{(q^6;q^6)_m}{(q^2;q^2)_{m+1}}.
	\]
	This proves the result.
\end{proof}

\begin{cor} \label{C=3}
	For $C=3$, the generating function for partitions in
	$\mathcal{P}_3$ is
	\[
	\sum_{\lambda\in\mathcal{P}_3} q^{|\lambda|}
	=
	\sum_{m\geq 0}
	\frac{
		q^{\frac{3m^2+7m+2}{2}}
	}{
		1-q^{2m+1}
	}
	\prod_{j=1}^{m}
	\left(
	1+q^{2j-1}
	+q^{2(2j-1)}
	+q^{3(2j-1)}
	\right).
	\]
	Equivalently,
	\[
	\sum_{\lambda\in\mathcal{P}_3} q^{|\lambda|}
	=
	\sum_{m\geq 0}
	q^{\frac{3m^2+7m+2}{2}}
	\frac{(q^4;q^8)_m}{(q;q^2)_{m+1}}.
	\]
\end{cor}

\begin{proof}
	Since $C=3$ is odd, there are no partitions in
	$\mathcal{P}_3$ with an even number of parts. Hence,
	\[
	\mathcal{P}_3=\mathcal{P}_3^{\mathrm{odd}}.
	\]
	
	Substituting $C=3$ into the generating function for
	$\mathcal{P}_C^{\mathrm{odd}}$ gives
	\[
	\sum_{\lambda\in\mathcal{P}_3} q^{|\lambda|}
	=
	\sum_{m\geq 0}
	\frac{
		q^{\frac{3m^2+7m+2}{2}}
	}{
		1-q^{2m+1}
	}
	\prod_{j=1}^{m}
	\left(
	\sum_{r=0}^{3}q^{r(2j-1)}
	\right).
	\]
	Since
	\[
	\sum_{r=0}^{3}q^{r(2j-1)}
	=
	1+q^{2j-1}
	+q^{2(2j-1)}
	+q^{3(2j-1)},
	\]
	the first formula follows.
	
	For the equivalent form, we use
	\[
	1+x+x^2+x^3=\frac{1-x^4}{1-x}.
	\]
	Taking $x=q^{2j-1}$ yields
	\[
	1+q^{2j-1}
	+q^{2(2j-1)}
	+q^{3(2j-1)}
	=
	\frac{1-q^{4(2j-1)}}{1-q^{2j-1}}.
	\]
	Therefore,
	\[
	\prod_{j=1}^{m}
	\left(
	1+q^{2j-1}
	+q^{2(2j-1)}
	+q^{3(2j-1)}
	\right)
	=
	\frac{
		\prod_{j=1}^{m}(1-q^{8j-4})
	}{
		\prod_{j=1}^{m}(1-q^{2j-1})
	}
	=
	\frac{(q^4;q^8)_m}{(q;q^2)_m}.
	\]
	Finally,
	\[
	(1-q^{2m+1})(q;q^2)_m=(q;q^2)_{m+1}.
	\]
	Substituting these identities into the first expression gives
	\[
	\sum_{\lambda\in\mathcal{P}_3} q^{|\lambda|}
	=
	\sum_{m\geq 0}
	q^{\frac{3m^2+7m+2}{2}}
	\frac{(q^4;q^8)_m}{(q;q^2)_{m+1}},
	\]
	as required.
\end{proof}

\begin{cor} \label{C=4}
	For $C=4$, the generating function for partitions in $\mathcal{P}_4$
	is
	\[
	\begin{aligned}
		\sum_{\lambda\in\mathcal{P}_4}q^{|\lambda|}
		={}&
		\sum_{m\geq0}
		\frac{q^{2m^2+4m+1}}{1-q^{2m+1}}
		\prod_{j=1}^{m}
		\left(
		1+q^{2j-1}+q^{2(2j-1)}
		+q^{3(2j-1)}+q^{4(2j-1)}
		\right)
		\\
		&+
		\sum_{m\geq0}
		\frac{q^{2(m+1)(m+2)}}{1-q^{2m+2}}
		\prod_{j=1}^{m}
		\left(
		1+q^{2j}+q^{4j}+q^{6j}+q^{8j}
		\right).
	\end{aligned}
	\]
	Equivalently,
	\[
	\begin{aligned}
		\sum_{\lambda\in\mathcal{P}_4}q^{|\lambda|}
		={}&
		\sum_{m\geq0}
		q^{2m^2+4m+1}
		\frac{(q^5;q^{10})_m}{(q;q^2)_{m+1}}
		+
		\sum_{m\geq0}
		q^{2(m+1)(m+2)}
		\frac{(q^{10};q^{10})_m}{(q^2;q^2)_{m+1}}.
	\end{aligned}
	\]
\end{cor}

\begin{proof}
	Since $C=4$ is even, both the odd- and even-length cases contribute.
	Substituting $C=4$ into Theorem~\ref{thm:GF} gives
	\[
	\sum_{\lambda\in\mathcal{P}_4^{\mathrm{odd}}}q^{|\lambda|}
	=
	\sum_{m\geq0}
	\frac{q^{2m^2+4m+1}}{1-q^{2m+1}}
	\prod_{j=1}^{m}
	\left(
	\sum_{r=0}^{4}q^{r(2j-1)}
	\right)
	\]
	and
	\[
	\sum_{\lambda\in\mathcal{P}_4^{\mathrm{even}}}q^{|\lambda|}
	=
	\sum_{m\geq0}
	\frac{q^{2(m+1)(m+2)}}{1-q^{2m+2}}
	\prod_{j=1}^{m}
	\left(
	\sum_{r=0}^{4}q^{2rj}
	\right).
	\]
	Since
	\[
	\sum_{r=0}^{4}q^{r(2j-1)}
	=
	1+q^{2j-1}+q^{2(2j-1)}
	+q^{3(2j-1)}+q^{4(2j-1)}
	\]
	and
	\[
	\sum_{r=0}^{4}q^{2rj}
	=
	1+q^{2j}+q^{4j}+q^{6j}+q^{8j},
	\]
	adding the odd- and even-length generating functions gives the first
	formula.
	
	For the equivalent form, we use
	\[
	1+x+x^2+x^3+x^4=\frac{1-x^5}{1-x}.
	\]
	Therefore,
	\[
	\prod_{j=1}^{m}
	\left(
	1+q^{2j-1}+q^{2(2j-1)}
	+q^{3(2j-1)}+q^{4(2j-1)}
	\right)
	=
	\frac{(q^5;q^{10})_m}{(q;q^2)_m},
	\]
	while
	\[
	\prod_{j=1}^{m}
	\left(
	1+q^{2j}+q^{4j}+q^{6j}+q^{8j}
	\right)
	=
	\frac{(q^{10};q^{10})_m}{(q^2;q^2)_m}.
	\]
	Combining these expressions with the factors
	\[
	\frac{1}{1-q^{2m+1}}
	\qquad\text{and}\qquad
	\frac{1}{1-q^{2m+2}},
	\]
	respectively, yields the stated $q$-Pochhammer forms.
\end{proof}

The special cases presented above reveal a common combinatorial
pattern. In each instance, the product part of the generating
function is precisely the generating function for a family of
partitions with bounded multiplicities. This observation suggests
that the generating functions obtained in Theorem~\ref{thm:GF} admit
a direct bijective interpretation. We now establish such bijections
for arbitrary values of the parameter $C$.

\section{Bijections with Bounded Multiplicity Partitions}

In this section, we establish explicit bijections between partitions
in $\mathcal{P}_C$ and suitable families of bounded multiplicity
partitions. These bijections provide a direct combinatorial
interpretation of the generating functions obtained in the previous
section.

For each integer $m\ge0$, let
$\mathcal{B}_{C,m}^{\mathrm{odd}}$
denote the set of partitions into odd parts not exceeding $2m+1$
such that each of the parts
\[
1,3,\ldots,2m-1
\]
occurs at most $C$ times, whereas the part $2m+1$ may occur with
arbitrary multiplicity.

For each integer $m\ge0$, let
\[
\mathcal{P}_C^{\mathrm{odd}}(2m+1)
=
\{\lambda\in\mathcal{P}_C:\ell(\lambda)=2m+1\},
\]
and
\[
\mathcal{P}_C^{\mathrm{even}}(2m)
=
\{\lambda\in\mathcal{P}_C:\ell(\lambda)=2m\}.
\]

\begin{theorem}\label{thm:odd-bijection}
	Let $C$ be a fixed nonnegative integer. Then, for every integer
	$m\ge0$, there exists an explicit bijection
	\[
	\Phi_{C,m}^{\mathrm{odd}}
	:
	\mathcal{P}_C^{\mathrm{odd}}(2m+1)
	\longrightarrow
	\mathcal{B}_{C,m}^{\mathrm{odd}}
	\]
	such that, if
	\[
	\mu=\Phi_{C,m}^{\mathrm{odd}}(\lambda),
	\]
	then
	\[
	|\lambda|
	=
	|\mu|
	+
	2m+1
	+
	C\binom{m+1}{2}.
	\]
\end{theorem}

\begin{proof}
	Let
	\[
	\lambda=(\lambda_1,\lambda_2,\ldots,\lambda_{2m+1})
	\in \mathcal{P}_C^{\mathrm{odd}}(2m+1),
	\]
	and let
	\[
	d_i=\lambda_{i+1}-\lambda_i,
	\qquad
	1\leq i\leq 2m,
	\]
	be its gap sequence. Since $\lambda\in\mathcal{P}_C$, we have
	\[
	d_i+d_{2m+1-i}=C,
	\qquad
	1\leq i\leq m.
	\]
	In particular,
	\[
	0\leq d_i\leq C,
	\qquad
	1\leq i\leq m.
	\]
	
	We define
	\[
	\mu=\Phi_{C,m}^{\mathrm{odd}}(\lambda)
	\]
	by prescribing the multiplicities of its parts as follows:
	\[
	\nu_{2j-1}(\mu)=d_{m+1-j},
	\qquad
	1\leq j\leq m,
	\]
	and
	\[
	\nu_{2m+1}(\mu)=\lambda_1-1,
	\]
	where $\nu_r(\mu)$ denotes the multiplicity of the part $r$ in $\mu$.
	Equivalently,
	\[
	\mu
	=
	1^{d_m}3^{d_{m-1}}\cdots
	(2m-1)^{d_1}(2m+1)^{\lambda_1-1}.
	\]
	
	Since
	\[
	0\leq d_i\leq C,
	\]
	each of the parts
	\[
	1,3,\ldots,2m-1
	\]
	occurs at most $C$ times in $\mu$. Moreover, the multiplicity of the
	part $2m+1$ is determined by the smallest part of $\lambda$, namely
	$
	\nu_{2m+1}(\mu)=\lambda_1-1.
	$
	Hence, 
	$
	\mu\in\mathcal{B}_{C,m}^{\mathrm{odd}},
	$
	so the map is well defined.
	
	We now construct the inverse map. Let
	\[
	\mu\in\mathcal{B}_{C,m}^{\mathrm{odd}}.
	\]
	Set
	\[
	\lambda_1=\nu_{2m+1}(\mu)+1
	\]
	and define
	\[
	d_i=\nu_{2m+1-2i}(\mu),
	\qquad
	1\leq i\leq m.
	\]
	Since each of the parts
	\[
	1,3,\ldots,2m-1
	\]
	occurs at most $C$ times in $\mu$, we have
	\[
	0\leq d_i\leq C.
	\]
	Define the remaining gaps by
	\[
	d_{2m+1-i}=C-d_i,
	\qquad
	1\leq i\leq m.
	\]
	Then all the gaps are nonnegative and satisfy
	\[
	d_i+d_{2m+1-i}=C.
	\]
	
	Starting from $\lambda_1$, define successively
	\[
	\lambda_{r+1}=\lambda_r+d_r,
	\qquad
	1\leq r\leq 2m.
	\]
	Since every $d_r$ is nonnegative, this produces a partition
	\[
	\lambda=(\lambda_1,\lambda_2,\ldots,\lambda_{2m+1})
	\]
	with exactly $2m+1$ positive parts. Moreover, its gap sequence
	satisfies the constant-sum condition, and therefore
	\[
	\lambda\in\mathcal{P}_C^{\mathrm{odd}}(2m+1).
	\]
	The two constructions recover the same smallest part and the same
	gap sequence, and hence are inverse to one another. Thus, 
	$
	\Phi_{C,m}^{\mathrm{odd}}
	$
	is a bijection.
	
	It remains to verify the relation between the weights. As established
	from the gap parametrization,
	\[
	|\lambda|
	=
	(2m+1)\lambda_1
	+
	C\binom{m+1}{2}
	+
	\sum_{j=1}^{m}(2j-1)d_{m+1-j}.
	\]
	On the other hand, by the definition of $\mu$,
	\[
	|\mu|
	=
	(2m+1)(\lambda_1-1)
	+
	\sum_{j=1}^{m}(2j-1)d_{m+1-j}.
	\]
	Subtracting these two identities gives
	\[
	|\lambda|-|\mu|
	=
	2m+1+C\binom{m+1}{2}.
	\]
	Therefore,
	\[
	|\lambda|
	=
	|\mu|
	+
	2m+1
	+
	C\binom{m+1}{2},
	\]
	as required.
\end{proof}

\begin{ex}
	Let $C=3$ and $m=2$, and consider the partition
	\[
	\lambda=(2,3,6,6,8).
	\]
	Its gap sequence is
	\[
	(d_1,d_2,d_3,d_4)=(1,3,0,2).
	\]
	Since
	\[
	d_1+d_4=1+2=3
	\]
	and
	\[
	d_2+d_3=3+0=3,
	\]
	we have
	\[
	\lambda\in\mathcal{P}_3^{\mathrm{odd}}(5).
	\]
	
	By the definition of the bijection
	$
	\Phi_{3,2}^{\mathrm{odd}},
	$
	the multiplicities of the parts of
	\[
	\mu=\Phi_{3,2}^{\mathrm{odd}}(\lambda)
	\]
	are given by
	\[
	\nu_1(\mu)=d_2=3,
	\qquad
	\nu_3(\mu)=d_1=1,
	\qquad
	\nu_5(\mu)=\lambda_1-1=1.
	\]
	Therefore,
	\[
	\mu=(1,1,1,3,5).
	\]
Notice that the parts $1$ and $3$ occur at most three times, as
	required by membership in $\mathcal{B}_{3,2}^{\mathrm{odd}}$, while the
	multiplicity of the part $5$ is not restricted by this membership
	condition; here it is determined by $\lambda_1$, namely
	$\nu_5(\mu)=\lambda_1-1=1$. Hence,
	$
	\mu\in\mathcal{B}_{3,2}^{\mathrm{odd}}.
	$
	
	Conversely, starting with
	\[
	\mu=(1,1,1,3,5),
	\]
	we recover
	\[
	\lambda_1=\nu_5(\mu)+1=2,
	\]
	and
	\[
	d_1=\nu_3(\mu)=1,
	\qquad
	d_2=\nu_1(\mu)=3.
	\]
	The remaining gaps are determined by the constant-sum condition:
	\[
	d_3=3-d_2=0,
	\qquad
	d_4=3-d_1=2.
	\]
	Thus,
	\[
	\lambda_2=\lambda_1+d_1=3,
	\]
	\[
	\lambda_3=\lambda_2+d_2=6,
	\]
	\[
	\lambda_4=\lambda_3+d_3=6,
	\]
	and
	\[
	\lambda_5=\lambda_4+d_4=8.
	\]
	Therefore, the inverse construction recovers
	\[
	\lambda=(2,3,6,6,8).
	\]
	
	Finally,
	\[
	|\lambda|=2+3+6+6+8=25
	\]
	and
	\[
	|\mu|=1+1+1+3+5=11.
	\]
	The weight relation in Theorem~\ref{thm:odd-bijection} is verified,
	since
	\[
	|\mu|+2m+1+C\binom{m+1}{2}
	=
	11+5+3\binom{3}{2}
	=
	11+5+9
	=
	25
	=
	|\lambda|.
	\]
\end{ex}

For each integer $m\geq1$, let
$\mathcal{B}_{C,m}^{\mathrm{even}}$
denote the set of partitions into even parts not exceeding $2m$
such that each of the parts
\[
2,4,\ldots,2m-2
\]
occurs at most $C$ times, whereas the part $2m$ may occur with
arbitrary multiplicity.

\noindent\textit{Remark.} 
	In Theorem~\ref{thm:GF}, partitions of even length
	were indexed as $\mathcal{P}_C^{\mathrm{even}}(2m+2)$, $m \geq 0$, so
	that the generating function could be summed starting from $m = 0$.
	For the bijection below it is more convenient to index a partition of
	even length directly by its length, so we write
	$\mathcal{P}_C^{\mathrm{even}}(2m)$ with $m \geq 1$ instead; the two
	conventions describe the same set of partitions, with $m$ shifted by
	one between them.

\begin{theorem}\label{thm:even-bijection}
	Let $C$ be a fixed nonnegative even integer. Then, for every integer
	$m\geq1$, there exists an explicit bijection
	\[
	\Phi_{C,m}^{\mathrm{even}}
	:
	\mathcal{P}_C^{\mathrm{even}}(2m)
	\longrightarrow
	\mathcal{B}_{C,m}^{\mathrm{even}}
	\]
	such that, if
	\[
	\mu=\Phi_{C,m}^{\mathrm{even}}(\lambda),
	\]
	then
	\[
	|\lambda|
	=
	|\mu|
	+
	2m
	+
	\frac{Cm^2}{2}.
	\]
\end{theorem}

\begin{proof}
	Let
	\[
	\lambda=(\lambda_1,\lambda_2,\ldots,\lambda_{2m})
	\in\mathcal{P}_C^{\mathrm{even}}(2m),
	\]
	and let
	\[
	d_i=\lambda_{i+1}-\lambda_i,
	\qquad
	1\leq i\leq 2m-1,
	\]
	be its gap sequence. Since $\lambda\in\mathcal{P}_C$, we have
	\[
	d_i+d_{2m-i}=C,
	\qquad
	1\leq i\leq m.
	\]
	In particular, for the middle gap,
	\[
	2d_m=C.
	\]
	Since $C$ is even, it follows that
	\[
	d_m=\frac{C}{2}.
	\]
	
	The gaps
	\[
	d_1,d_2,\ldots,d_{m-1}
	\]
	may be chosen freely from $\{0,1,\ldots,C\}$, while the remaining
	gaps are uniquely determined by
	\[
	d_{2m-i}=C-d_i,
	\qquad
	1\leq i\leq m-1.
	\]
	
	We define
	\[
	\mu=\Phi_{C,m}^{\mathrm{even}}(\lambda)
	\]
	by prescribing the multiplicities of its parts as follows:
	\[
	\nu_{2j}(\mu)=d_{m-j},
	\qquad
	1\leq j\leq m-1,
	\]
	and
	\[
	\nu_{2m}(\mu)=\lambda_1-1,
	\]
	where $\nu_r(\mu)$ denotes the multiplicity of the part $r$ in $\mu$.
	Equivalently,
	\[
	\mu
	=
	2^{d_{m-1}}4^{d_{m-2}}\cdots
	(2m-2)^{d_1}(2m)^{\lambda_1-1}.
	\]
	
	Since
	\[
	0\leq d_i\leq C,
	\qquad
	1\leq i\leq m-1,
	\]
	each of the parts
	\[
	2,4,\ldots,2m-2
	\]
	occurs at most $C$ times in $\mu$. Moreover, the multiplicity of the part $2m$ is determined by the smallest part of $\lambda$, namely 	$\nu_{2m}(\mu)=\lambda_1-1$.	 Hence,
	$
	\mu\in\mathcal{B}_{C,m}^{\mathrm{even}},
	$
	and therefore the map is well defined.
	
	We now construct the inverse map. Let
	$
	\mu\in\mathcal{B}_{C,m}^{\mathrm{even}}
	$.
	Set
	\[
	\lambda_1=\nu_{2m}(\mu)+1
	\]
	and define
	\[
	d_i=\nu_{2m-2i}(\mu),
	\qquad
	1\leq i\leq m-1.
	\]
	Since each of the parts
	\[
	2,4,\ldots,2m-2
	\]
	occurs at most $C$ times in $\mu$, we have
	\[
	0\leq d_i\leq C.
	\]
	Set
	\[
	d_m=\frac{C}{2},
	\]
	and define the remaining gaps by
	\[
	d_{2m-i}=C-d_i,
	\qquad
	1\leq i\leq m-1.
	\]
	Then all the gaps are nonnegative and satisfy
	\[
	d_i+d_{2m-i}=C,
	\qquad
	1\leq i\leq m.
	\]
	
	Starting from $\lambda_1$, define successively
	\[
	\lambda_{r+1}=\lambda_r+d_r,
	\qquad
	1\leq r\leq 2m-1.
	\]
	Since every $d_r$ is nonnegative, this produces a partition
	\[
	\lambda=(\lambda_1,\lambda_2,\ldots,\lambda_{2m})
	\]
	with exactly $2m$ positive parts. Its gap sequence satisfies the
	constant-sum condition, and hence
	\[
	\lambda\in\mathcal{P}_C^{\mathrm{even}}(2m).
	\]
	
	The two constructions recover the same smallest part and the same
	gap sequence, and therefore they are inverse to one another. Thus,
	$
	\Phi_{C,m}^{\mathrm{even}}
	$
	is a bijection.
	
	It remains to verify the relation between the weights. From the gap
	parametrization, we have
	\[
	|\lambda|
	=
	2m\lambda_1
	+
	\frac{Cm^2}{2}
	+
	\sum_{j=1}^{m-1}2j\,d_{m-j}.
	\]
	On the other hand, by the definition of $\mu$,
	\[
	|\mu|
	=
	2m(\lambda_1-1)
	+
	\sum_{j=1}^{m-1}2j\,d_{m-j}.
	\]
	Subtracting the second identity from the first gives
	\[
	|\lambda|-|\mu|
	=
	2m+\frac{Cm^2}{2}.
	\]
	Therefore,
	\[
	|\lambda|
	=
	|\mu|
	+
	2m
	+
	\frac{Cm^2}{2},
	\]
	as required.
\end{proof}

\begin{ex}
	Let $C=4$ and $m=3$, and consider the partition
	\[
	\lambda=(2,3,7,9,9,12).
	\]
	Its gap sequence is
	\[
	(d_1,d_2,d_3,d_4,d_5)=(1,4,2,0,3).
	\]
	We have
	\[
	d_1+d_5=1+3=4,
	\]
	\[
	d_2+d_4=4+0=4,
	\]
	and
	\[
	2d_3=2\cdot2=4.
	\]
	Therefore,
	\[
	\lambda\in\mathcal{P}_4^{\mathrm{even}}(6).
	\]
	
	By the definition of the bijection
	$
	\Phi_{4,3}^{\mathrm{even}},
	$
	the multiplicities of the parts of
	\[
	\mu=\Phi_{4,3}^{\mathrm{even}}(\lambda)
	\]
	are given by
	\[
	\nu_2(\mu)=d_2=4,
	\qquad
	\nu_4(\mu)=d_1=1,
	\qquad
	\nu_6(\mu)=\lambda_1-1=1.
	\]
	Thus,
	\[
	\mu=(2,2,2,2,4,6).
	\]
	The parts $2$ and $4$ occur at most four times. Moreover, the
	multiplicity of the part $6$ is determined by the smallest part of
	$\lambda$, namely
	$
	\nu_6(\mu)=\lambda_1-1=1.
	$ Hence,
	$
	\mu\in\mathcal{B}_{4,3}^{\mathrm{even}}.
	$
	
	Conversely, starting with
	\[
	\mu=(2,2,2,2,4,6),
	\]
	we recover the smallest part as
	\[
	\lambda_1=\nu_6(\mu)+1=2.
	\]
	The first two gaps are
	\[
	d_1=\nu_4(\mu)=1,
	\qquad
	d_2=\nu_2(\mu)=4.
	\]
	The middle gap is fixed by
	\[
	d_3=\frac{C}{2}=2.
	\]
	The remaining gaps are determined by the constant-sum condition:
	\[
	d_4=C-d_2=4-4=0,
	\]
	and
	\[
	d_5=C-d_1=4-1=3.
	\]
	Starting from $\lambda_1=2$, we obtain
	\[
	\lambda_2=\lambda_1+d_1=3,
	\]
	\[
	\lambda_3=\lambda_2+d_2=7,
	\]
	\[
	\lambda_4=\lambda_3+d_3=9,
	\]
	\[
	\lambda_5=\lambda_4+d_4=9,
	\]
	and
	\[
	\lambda_6=\lambda_5+d_5=12.
	\]
	Thus, the inverse construction recovers
	\[
	\lambda=(2,3,7,9,9,12).
	\]
	
	Finally,
	\[
	|\lambda|=2+3+7+9+9+12=42
	\]
	and
	\[
	|\mu|=2+2+2+2+4+6=18.
	\]
	The weight relation in Theorem~\ref{thm:even-bijection} is verified,
	since
	\[
	\begin{aligned}
		|\mu|+2m+\frac{Cm^2}{2}
		&=
		18+6+\frac{4\cdot 3^2}{2}\\
		&=
		18+6+18\\
		&=
		42\\
		&=
		|\lambda|.
	\end{aligned}
	\]
\end{ex}

\section{Ferrers Diagram Characterizations}

In this section, we establish geometric characterizations of the
partitions in $\mathcal{P}_C$ by means of their Ferrers diagrams.
These characterizations reveal that the constant-sum gap condition is
equivalent to a natural symmetry property of suitably modified Ferrers
diagrams.

\begin{theorem}[Ferrers diagram characterization, odd length] \label{thm:ferrers-odd}
	Let
	\[
	\lambda=(\lambda_1,\lambda_2,\ldots,\lambda_{2m+1})
	\]
	be an integer partition, and let $F(\lambda)$ denote its Ferrers
	diagram.
	
	Then
	\[
	\lambda\in\mathcal{P}_C^{\mathrm{odd}}(2m+1)
	\]
	if and only if the Ferrers diagram $F(\lambda)$ has the following
	property:
	
	For each $r=1,\ldots,m$, remove the rightmost $rC$ cells from the
	$(m+1+r)$-th row. The resulting diagram has row lengths
	\[
	\lambda_1,\lambda_2,\ldots,\lambda_m,\lambda_{m+1},
	\lambda_m,\ldots,\lambda_2,\lambda_1.
	\]
	Equivalently, after these removals, the Ferrers diagram becomes
	symmetric with respect to its middle row.
\end{theorem}

\begin{proof}
	Let
	\[
	d_i=\lambda_{i+1}-\lambda_i,
	\qquad
	1\leq i\leq 2m,
	\]
	be the gap sequence of $\lambda$.
	
	Suppose first that
	\[
	\lambda\in\mathcal{P}_C^{\mathrm{odd}}(2m+1).
	\]
	Then
	\[
	d_i+d_{2m+1-i}=C,
	\qquad
	1\leq i\leq m.
	\]
	
	Fix $r\in\{1,\ldots,m\}$. The difference between the lengths of the
	two rows situated symmetrically with respect to the middle row is
	\[
	\lambda_{m+1+r}-\lambda_{m+1-r}.
	\]
	Using the gap sequence, we have
	\[
	\lambda_{m+1+r}-\lambda_{m+1-r}
	=
	\sum_{j=m+1-r}^{m+r}d_j.
	\]
	The terms in this sum may be paired symmetrically as
	\[
	(d_{m+1-r},d_{m+r}),
	\,
	(d_{m+2-r},d_{m+r-1}),
	\ldots,
	(d_m,d_{m+1}).
	\]
	Each pair has sum $C$, since
	\[
	d_j+d_{2m+1-j}=C.
	\]
	There are exactly $r$ such pairs. Therefore,
	\[
	\lambda_{m+1+r}-\lambda_{m+1-r}=rC.
	\]
	Hence, after removing the rightmost $rC$ cells from the
	$(m+1+r)$-th row, its new length is
	\[
	\lambda_{m+1+r}-rC
	=
	\lambda_{m+1-r}.
	\]
	
	Since this holds for every $r=1,\ldots,m$, the row lengths of the
	modified diagram are
	\[
	\lambda_1,\lambda_2,\ldots,\lambda_m,\lambda_{m+1},
	\lambda_m,\ldots,\lambda_2,\lambda_1.
	\]
	Thus, the resulting diagram is symmetric with respect to its middle
	row.
	
	Conversely, suppose that, for each $r=1,\ldots,m$, removing the
	rightmost $rC$ cells from the $(m+1+r)$-th row produces a diagram
	whose row lengths are
	\[
	\lambda_1,\lambda_2,\ldots,\lambda_m,\lambda_{m+1},
	\lambda_m,\ldots,\lambda_2,\lambda_1.
	\]
	Then
	\[
	\lambda_{m+1+r}-rC
	=
	\lambda_{m+1-r},
	\qquad
	1\leq r\leq m,
	\]
	or equivalently,
	\[
	\lambda_{m+1+r}-\lambda_{m+1-r}=rC.
	\]
	
	For $1\leq r\leq m$, subtracting the identity corresponding to
	$r-1$ from the one corresponding to $r$, with the case $r=1$
	interpreted using the trivial identity
	\[
	\lambda_{m+1}-\lambda_{m+1}=0,
	\]
	gives
	\[
	\begin{aligned}
		C
		&=
		\bigl(\lambda_{m+1+r}-\lambda_{m+1-r}\bigr)
		-
		\bigl(\lambda_{m+r}-\lambda_{m+2-r}\bigr)\\
		&=
		\bigl(\lambda_{m+1+r}-\lambda_{m+r}\bigr)
		+
		\bigl(\lambda_{m+2-r}-\lambda_{m+1-r}\bigr)\\
		&=
		d_{m+r}+d_{m+1-r}.
	\end{aligned}
	\]
	Setting
	\[
	i=m+1-r,
	\]
	we obtain
	\[
	d_i+d_{2m+1-i}=C,
	\qquad
	1\leq i\leq m.
	\]
	Therefore,
	\[
	\lambda\in\mathcal{P}_C^{\mathrm{odd}}(2m+1).
	\]
	This completes the proof.
\end{proof}

\begin{ex}
	Let $C=3$ and $m=2$, and consider the partition
	\[
	\lambda=(2,3,6,6,8).
	\]
	Its gap sequence is
	\[
	(d_1,d_2,d_3,d_4)=(1,3,0,2).
	\]
	Since
	\[
	d_1+d_4=1+2=3
	\]
	and
	\[
	d_2+d_3=3+0=3,
	\]
	we have
	\[
	\lambda\in\mathcal{P}_3^{\mathrm{odd}}(5).
	\]
	
	We now apply the Ferrers-diagram modification described in
	Theorem~\ref{thm:ferrers-odd}. Since $m=2$, we remove
	\[
	rC
	\]
	cells from the $(m+1+r)$-th row for $r=1,2$.
	
	For $r=1$, we remove
	\[
	1\cdot C=3
	\]
	cells from the fourth row. Its new length is
	\[
	\lambda_4-3=6-3=3=\lambda_2.
	\]
	
	For $r=2$, we remove
	\[
	2\cdot C=6
	\]
	cells from the fifth row. Its new length is
	\[
	\lambda_5-6=8-6=2=\lambda_1.
	\]
	
	Thus, the row lengths of the modified diagram are
	\[
	(2,3,6,3,2),
	\]
	which are symmetric with respect to the middle row. Hence, the
	Ferrers diagram of $\lambda$ satisfies the geometric condition in
	Theorem~\ref{thm:ferrers-odd}.
	
	The figure \ref{oddlengthexample} provides a graphical illustration of the above
	example, showing the Ferrers diagram of the partition together with
	the cell deletions described in
	Theorem~\ref{thm:ferrers-odd}.

\begin{figure} 
	\begin{tikzpicture}[scale=0.50, every node/.style={font=\small}]
		
		\foreach \x in {0,...,1} { \draw (\x,0) rectangle (\x+1,1); }
		\foreach \x in {0,...,2} { \draw (\x,1) rectangle (\x+1,2); }
		\foreach \x in {0,...,5} { \draw (\x,2) rectangle (\x+1,3); }
		\foreach \x in {0,...,2} { \draw (\x,3) rectangle (\x+1,4); }
		\foreach \x in {3,...,5} { \draw[fill=gray!30] (\x,3) rectangle (\x+1,4); }
		\foreach \x in {0,...,1} { \draw (\x,4) rectangle (\x+1,5); }
		\foreach \x in {2,...,7} { \draw[fill=gray!30] (\x,4) rectangle (\x+1,5); }
		
		\node[anchor=east] at (-0.3,0.5) {$\lambda_1=2$};
		\node[anchor=east] at (-0.3,1.5) {$\lambda_2=3$};
		\node[anchor=east] at (-0.3,2.5) {$\lambda_3=6$};
		\node[anchor=east] at (-0.3,3.5) {$\lambda_4=6$};
		\node[anchor=east] at (-0.3,4.5) {$\lambda_5=8$};
		
		\node at (4,6) {$F(\lambda)$};
		\node at (4,-1.3) {$(2,3,6,6,8)$};
		
		\draw[->] (5,4.5) -- (8.5,4.8);
		\node[anchor=west] at (8.5,4.8) {$6 = 2C$ cells};
		
		\draw[->] (4.5,3.5) -- (8.5,3.0);
		\node[anchor=west] at (8.5,3.0) {$3 = C$ cells};
		
		\node[anchor=west, align=left] at (8.5,1.3) {remove the\\ shaded cells};
		
		\draw[->, thick] (13,2.5) -- (13.8,2.5);
		
		\begin{scope}[xshift=15cm]
			\foreach \x in {0,...,1} { \draw (\x,0) rectangle (\x+1,1); }
			\foreach \x in {0,...,2} { \draw (\x,1) rectangle (\x+1,2); }
			\foreach \x in {0,...,5} { \draw (\x,2) rectangle (\x+1,3); }
			\foreach \x in {0,...,2} { \draw (\x,3) rectangle (\x+1,4); }
			\foreach \x in {0,...,1} { \draw (\x,4) rectangle (\x+1,5); }
			
			\node[anchor=east] at (-0.3,0.5) {$2$};
			\node[anchor=east] at (-0.3,1.5) {$3$};
			\node[anchor=east] at (-0.3,2.5) {$6$};
			\node[anchor=east] at (-0.3,3.5) {$3$};
			\node[anchor=east] at (-0.3,4.5) {$2$};
			
			\node at (3,6) {Modified diagram};
			\node at (3,-1.3) {$(2,3,6,3,2)$};
			
			\draw[dashed] (-0.3,2.5) -- (6.8,2.5);
			\node[anchor=west, align=left] at (7.0,2.5) {axis of symmetry\\ through the middle \\ row};
		\end{scope}
		
		\draw (5,-3.2) rectangle (5.7,-2.5);
		\node[anchor=west] at (5.7,-2.85) {retained cells};
		
		\draw[fill=gray!30] (10.9,-3.2) rectangle (11.6,-2.5);
		\node[anchor=west] at (11.5,-2.85) {removed cells};
		
	\end{tikzpicture}
    \caption{Ferrers-diagram modification for
	$\lambda=(2,3,6,6,8)\in\mathcal P_3^{\mathrm{odd}}(5)$.
	Removing the shaded cells transforms $F(\lambda)$ into a diagram
	symmetric with respect to its middle row.}

		\label{oddlengthexample}
\end{figure}
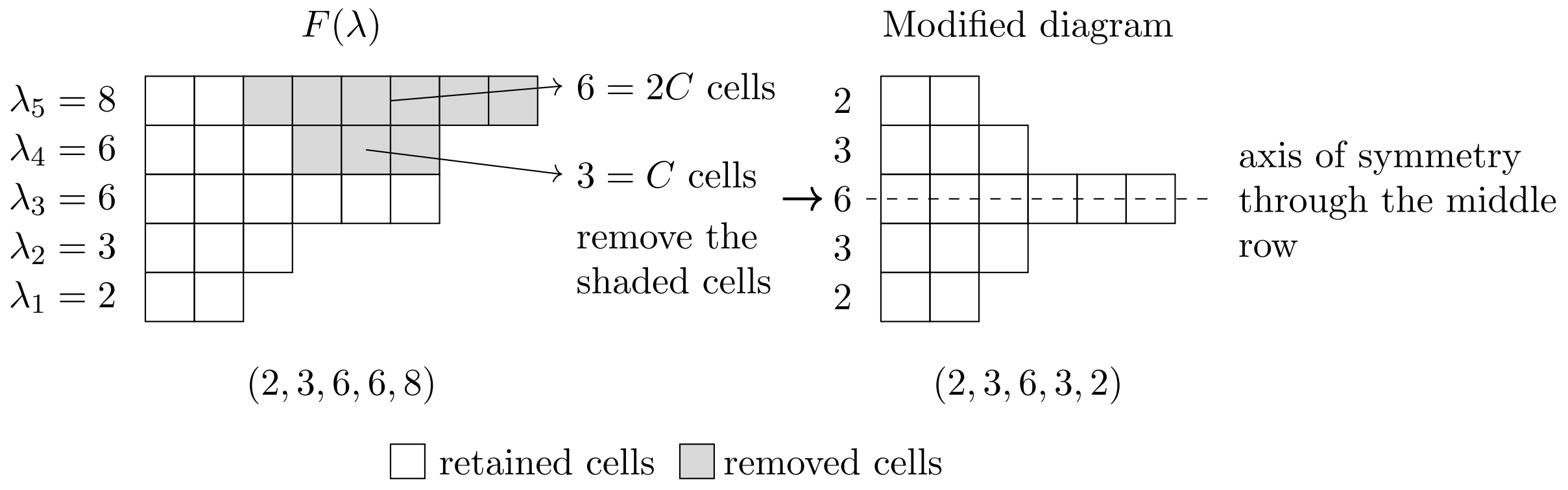

\end{ex}

\begin{theorem}[Ferrers diagram characterization, even length]
	\label{thm:ferrers-even}
	Let $C=2D$ be a fixed nonnegative even integer, and let
	\[
	\lambda=(\lambda_1,\lambda_2,\ldots,\lambda_{2m})
	\]
	be an integer partition. Let $F(\lambda)$ denote its Ferrers diagram.
	
	Then
	\[
	\lambda\in\mathcal{P}_C^{\mathrm{even}}(2m)
	\]
	if and only if the Ferrers diagram $F(\lambda)$ has the following
	property:
	
	For each $r=1,\ldots,m$, remove the rightmost
	\[
	(2r-1)D
	\]
	cells from the $(m+r)$-th row. The resulting diagram has row lengths
	\[
	\lambda_1,\lambda_2,\ldots,\lambda_{m-1},\lambda_m,
	\lambda_m,\lambda_{m-1},\ldots,\lambda_2,\lambda_1.
	\]
	Equivalently, after these removals, the resulting diagram is symmetric
	with respect to the horizontal line separating its two middle rows.
\end{theorem}

\begin{proof}
	Let
	\[
	d_i=\lambda_{i+1}-\lambda_i,
	\qquad
	1\leq i\leq 2m-1,
	\]
	be the gap sequence of $\lambda$.
	
	Suppose first that
	\[
	\lambda\in\mathcal{P}_C^{\mathrm{even}}(2m),
	\]
	where $C=2D$. Then
	\[
	d_i+d_{2m-i}=2D,
	\qquad
	1\leq i\leq m.
	\]
	In particular, taking $i=m$, we obtain
	\[
	2d_m=2D,
	\]
	and hence
	\[
	d_m=D.
	\]
	
	Fix $r\in\{1,\ldots,m\}$. The difference between the lengths of the
	two rows situated symmetrically with respect to the horizontal line
	separating the two middle rows is
	\[
	\lambda_{m+r}-\lambda_{m+1-r}.
	\]
	In terms of the gaps,
	\[
	\lambda_{m+r}-\lambda_{m+1-r}
	=
	\sum_{j=m+1-r}^{m+r-1}d_j.
	\]
	The terms in this sum consist of the middle gap $d_m$ together with
	the $r-1$ symmetric pairs
	\[
	(d_{m+1-r},d_{m+r-1}),
	\,
	(d_{m+2-r},d_{m+r-2}),
	\ldots,
	(d_{m-1},d_{m+1}).
	\]
	Each symmetric pair has sum $2D$, while $d_m=D$. Therefore,
	\[
	\begin{aligned}
		\lambda_{m+r}-\lambda_{m+1-r}
		&=(r-1)(2D)+D\\
		&=(2r-1)D.
	\end{aligned}
	\]
	It follows that, after removing the rightmost $(2r-1)D$ cells from
	the $(m+r)$-th row, its new length is
	\[
	\lambda_{m+r}-(2r-1)D
	=
	\lambda_{m+1-r}.
	\]
	Since this holds for every $r=1,\ldots,m$, the row lengths of the
	modified diagram are
	\[
	\lambda_1,\lambda_2,\ldots,\lambda_{m-1},\lambda_m,
	\lambda_m,\lambda_{m-1},\ldots,\lambda_2,\lambda_1.
	\]
	Thus, the resulting diagram is symmetric with respect to the
	horizontal line separating its two middle rows.
	
	Conversely, suppose that, for every $r=1,\ldots,m$, removing the
	rightmost $(2r-1)D$ cells from the $(m+r)$-th row produces a diagram
	with row lengths
	\[
	\lambda_1,\lambda_2,\ldots,\lambda_{m-1},\lambda_m,
	\lambda_m,\lambda_{m-1},\ldots,\lambda_2,\lambda_1.
	\]
	Then
	\[
	\lambda_{m+r}-(2r-1)D
	=
	\lambda_{m+1-r},
	\qquad
	1\leq r\leq m,
	\]
	or equivalently,
	\[
	\lambda_{m+r}-\lambda_{m+1-r}
	=
	(2r-1)D.
	\]
	
	For $r=1$, this gives
	\[
	\lambda_{m+1}-\lambda_m=D,
	\]
	and hence
	\[
	d_m=D.
	\]
	Therefore,
	\[
	2d_m=2D=C.
	\]
	
	Now let $2\leq r\leq m$. Subtracting the identity corresponding to
	$r-1$ from the identity corresponding to $r$, we obtain
	\[
	\begin{aligned}
		2D
		&=
		\bigl(\lambda_{m+r}-\lambda_{m+1-r}\bigr)
		-
		\bigl(\lambda_{m+r-1}-\lambda_{m+2-r}\bigr)\\
		&=
		\bigl(\lambda_{m+r}-\lambda_{m+r-1}\bigr)
		+
		\bigl(\lambda_{m+2-r}-\lambda_{m+1-r}\bigr)\\
		&=
		d_{m+r-1}+d_{m+1-r}.
	\end{aligned}
	\]
	Setting
	\[
	i=m+1-r,
	\]
	we obtain
	\[
	d_i+d_{2m-i}=2D=C,
	\qquad
	1\leq i\leq m-1.
	\]
	Together with
	\[
	2d_m=C,
	\]
	this shows that
	\[
	d_i+d_{2m-i}=C,
	\qquad
	1\leq i\leq m.
	\]
	Hence, 
	$
	\lambda\in\mathcal{P}_C^{\mathrm{even}}(2m)
	$.
	This completes the proof.
\end{proof}

\begin{ex}
	Let $C=4$, so that $D=C/2=2$, and let $m=3$. Consider the
	partition
	\[
	\lambda=(2,3,7,9,9,12).
	\]
	Its gap sequence is
	\[
	(d_1,d_2,d_3,d_4,d_5)=(1,4,2,0,3).
	\]
	Since
	\[
	d_1+d_5=1+3=4,
	\]
	\[
	d_2+d_4=4+0=4,
	\]
	and
	\[
	2d_3=2\cdot2=4,
	\]
	we have
	\[
	\lambda\in\mathcal{P}_4^{\mathrm{even}}(6).
	\]
	
	We now apply the Ferrers-diagram modification described in
	Theorem~\ref{thm:ferrers-even}. Since $m=3$, for each
	$r=1,2,3$ we remove the rightmost
	\[
	(2r-1)D
	\]
	cells from the $(m+r)$-th row.
	
	For $r=1$, we remove
	\[
	(2\cdot1-1)D=2
	\]
	cells from the fourth row. Its new length is
	\[
	\lambda_4-2=9-2=7=\lambda_3.
	\]
	
	For $r=2$, we remove
	\[
	(2\cdot2-1)D=6
	\]
	cells from the fifth row. Its new length is
	\[
	\lambda_5-6=9-6=3=\lambda_2.
	\]
	
	For $r=3$, we remove
	\[
	(2\cdot3-1)D=10
	\]
	cells from the sixth row. Its new length is
	\[
	\lambda_6-10=12-10=2=\lambda_1.
	\]
	
	Thus, the row lengths of the modified diagram are
	\[
	(2,3,7,7,3,2).
	\]
	These row lengths are symmetric with respect to the horizontal line
	separating the third and fourth rows. Hence, the Ferrers diagram of
	$\lambda$ satisfies the geometric condition in
	Theorem~\ref{thm:ferrers-even}.
	
The figure \ref{evenlengthexample}, using Ferrers diagrams, provides a visualization of this example.

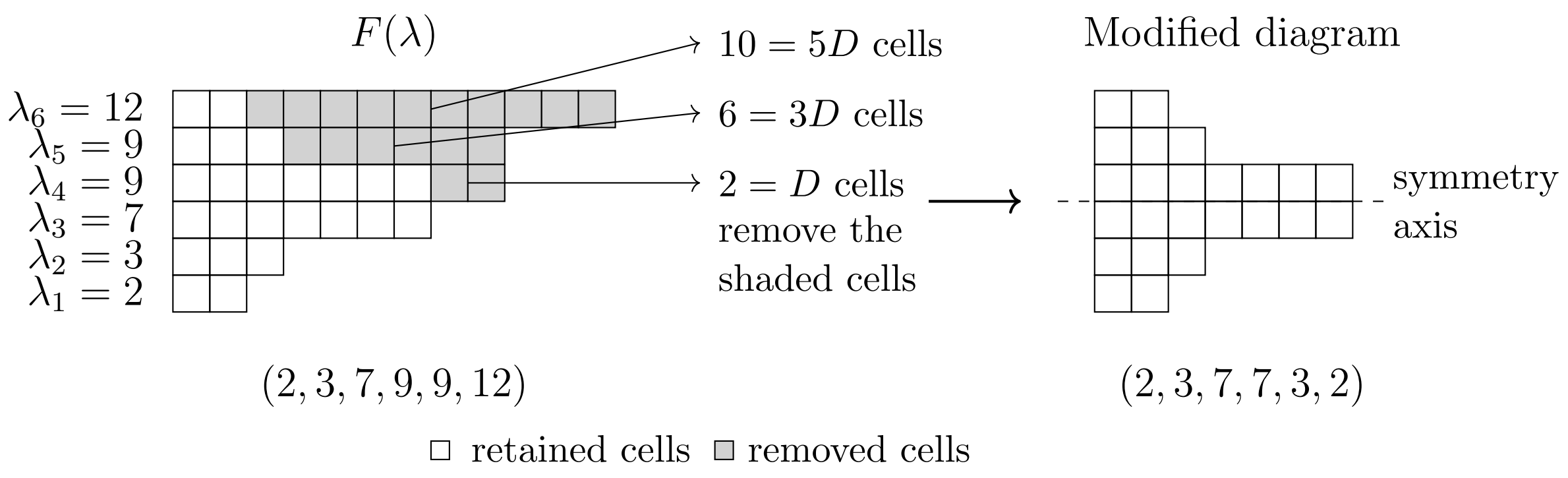
\begin{figure}[ht]
	\centering
	
	\begin{tikzpicture}[x=.37cm,y=.37cm]
		
		
		\node at (6,2.5) {$F(\lambda)$};
		
		\foreach \x in {0,1}
		\draw (\x,0) rectangle ++(1,1);
		\foreach \x in {2,...,11}
		\filldraw[fill=gray!35] (\x,0) rectangle ++(1,1);
		
		\foreach \x in {0,1,2}
		\draw (\x,-1) rectangle ++(1,1);
		\foreach \x in {3,...,8}
		\filldraw[fill=gray!35] (\x,-1) rectangle ++(1,1);
		
		\foreach \x in {0,...,6}
		\draw (\x,-2) rectangle ++(1,1);
		\foreach \x in {7,8}
		\filldraw[fill=gray!35] (\x,-2) rectangle ++(1,1);
		
		\foreach \x in {0,...,6}
		\draw (\x,-3) rectangle ++(1,1);
		
		\foreach \x in {0,...,2}
		\draw (\x,-4) rectangle ++(1,1);
		
		\foreach \x in {0,1}
		\draw (\x,-5) rectangle ++(1,1);
		
		
		\node[left] at (-.4,.5) {$\lambda_6=12$};
		\node[left] at (-.4,-.5) {$\lambda_5=9$};
		\node[left] at (-.4,-1.5) {$\lambda_4=9$};
		\node[left] at (-.4,-2.5) {$\lambda_3=7$};
		\node[left] at (-.4,-3.5) {$\lambda_2=3$};
		\node[left] at (-.4,-4.5) {$\lambda_1=2$};
		
		
		\draw[->] (7,.5) -- (14.3,2.3);
		\node[right] at (14.4,2.3) {\small $10=5D$ cells};
		
		\draw[->] (6,-.5) -- (14.3,.4);
		\node[right] at (14.4,.4) {\small $6=3D$ cells};
		
		\draw[->] (8,-1.5) -- (14.3,-1.5);
		\node[right] at (14.4,-1.5) {\small $2=D$ cells};
		
		\node[right, align=left] at (14.4,-3.4) {\small remove the\\ \small shaded cells};
		
		
		\draw[->,thick] (20.5,-2) -- (23,-2);
		
		
		\node at (29,2.5) {Modified diagram};
		
		\foreach \x in {0,1}
		\draw (25+\x,0) rectangle ++(1,1);
		
		\foreach \x in {0,1,2}
		\draw (25+\x,-1) rectangle ++(1,1);
		
		\foreach \x in {0,...,6}
		\draw (25+\x,-2) rectangle ++(1,1);
		
		\foreach \x in {0,...,6}
		\draw (25+\x,-3) rectangle ++(1,1);
		
		\foreach \x in {0,...,2}
		\draw (25+\x,-4) rectangle ++(1,1);
		
		\foreach \x in {0,1}
		\draw (25+\x,-5) rectangle ++(1,1);
		
		
		\draw[dashed] (24,-2) -- (33,-2);
		\node[right, align=left] at (32.7,-2) {\small symmetry\\ \small axis};
		
		
		\node at (6,-7) {$(2,3,7,9,9,12)$};
		\node at (29,-7) {$(2,3,7,7,3,2)$};
		
		
		\draw (7,-9) rectangle ++(.5,.5);
		\node[right] at (7.7,-8.70) {\small retained cells};
		
		\filldraw[fill=gray!35] (14.7,-9) rectangle ++(.5,.5);
		\node[right] at (15.2,-8.70) {\small removed cells};
		
	\end{tikzpicture}
	
	\caption{Ferrers-diagram modification for
		$\lambda=(2,3,7,9,9,12)\in\mathcal P^{\mathrm{even}}_4(6)$.
		Removing $5D$, $3D$, and $D$ cells from the first, second, and third
		rows, respectively, produces a diagram symmetric with respect to the
		horizontal line separating the two middle rows.}
	\label{evenlengthexample}
\end{figure}	

\end{ex}

\begin{cor}[Deleted area]
	Let $F(\lambda)$ be the Ferrers diagram of a partition
	\[
	\lambda\in\mathcal{P}_C.
	\]
	
	\begin{enumerate}[(i)]
		\item If
		\[
		\lambda\in\mathcal{P}_C^{\mathrm{odd}}(2m+1),
		\]
		then the total number of deleted cells in
		Theorem~\ref{thm:ferrers-odd} is
		\[
		C\binom{m+1}{2}.
		\]
		
		\item If
		\[
		\lambda\in\mathcal{P}_C^{\mathrm{even}}(2m),
		\]
		where $C$ is even, then the total number of deleted cells in
		Theorem~\ref{thm:ferrers-even} is
		\[
		\frac{Cm^{2}}{2}.
		\]
	\end{enumerate}
\end{cor}

\begin{proof}
	In the odd case, the deleted cells form strips of lengths
	\[
	C,2C,\ldots,mC,
	\]
	whose total area is
	\[
	C(1+2+\cdots+m)
	=
	C\binom{m+1}{2}.
	\]
	
	In the even case, writing $C=2D$, the deleted strips have lengths
	\[
	D,3D,\ldots,(2m-1)D.
	\]
	Since
	\[
	1+3+\cdots+(2m-1)=m^{2},
	\]
	their total area is
	\[
	Dm^{2}
	=
	\frac{Cm^{2}}{2}.
	\]
	This completes the proof.
\end{proof}

\section{Conclusion}

The results of this paper show that the constant-sum gap condition,
though defined by a single linear relation on the gap sequence, forces
a rigid combinatorial structure: the generating function splits
cleanly by parity of length, each parity class is in explicit bijection
with a family of bounded-multiplicity partitions, and this same
structure is visible geometrically as a controlled symmetrization of
the Ferrers diagram. The coincidence between the area removed in the
Ferrers-diagram description and the weight shift appearing in the
bijections is not accidental: both quantities arise from summing the
same $C$-dependent contributions over the gap sequence, so the two
results are, at bottom, two descriptions of the same parametrization.

Several questions remain open. The generating functions obtained here
do not appear to simplify beyond a $q$-Pochhammer form, and it would be
of interest to know whether a product formula or modular interpretation
exists for general $C$. Beyond this, throughout this paper $C$ has been
fixed, and it would be natural to ask what generating function results
from summing, or otherwise combining, the families $\mathcal{P}_C$ over
all $C \geq 0$. On a
different note, it would also be natural to investigate more closely
the asymmetry between the odd- and even-length cases -- most notably
the vanishing of $\mathcal{P}_C^{\mathrm{even}}$ for odd $C$. We leave these questions for future work.

\bibliographystyle{amsplain}

\end{document}